\documentclass[pssymbold,secthm]{elsart}
\usepackage[psamsfonts]{amsfonts}
\usepackage[psamsfonts]{amssymb}
\usepackage{mathrsfs}
\usepackage{amsmath}

\def\CD{Cartan-Dieudonn\'e }
\def\ort{orthogonal }

\begin{document}

\begin{frontmatter}

  \title{An algorithm for the Cartan-Dieudonn\'e theorem on
    generalized scalar product spaces}

  \author[IPN]{M.A. Rodr\'{\i}guez-Andrade}, \ead{marco@esfm.ipn.mx}
  \author[UAMA]{G. Arag\'on-Gonz\'alez}, \ead{gag@correo.azc.uam.mx}
  \author[UNAM]{J.L. Arag\'on}, and  \ead{aragon@fata.unam.mx}
  \author[UAM]{Luis Verde-Star} \ead{verde@xanum.uam.mx}

  \address[IPN]{Departamento de Matem\'aticas, Escuela Superior de
    F\'{\i}sica y Matem\'aticas, IPN, M\'exico D.F. 07300, M\'exico
    and Departamento de Matem\'atica Educativa Centro de
    Investigaci\'on y Estudios Avanzados, IPN, M\'exico D.F. 07360,
    M\'exico.}

  \address[UAMA]{Programa de Desarrollo Profesional en Automatizaci\'on,
  Universidad Aut\'onoma Metropolitana, Azcapotzalco, San Pablo 180,
  Colonia Reynosa-Tamaulipas, M\'exico D.F. 02200, M\'exico.}

 \address[UNAM]{Centro de F\'{\i}sica Aplicada y Tecnolog\'{\i}a
    Avanzada, Universidad Nacional Aut\'onoma de M\'exico, Apartado
    1-1010, Quer\'etaro 76000, M\'exico.}
  
   \address[UAM]{Departamento de Matem\'aticas, Universidad Aut\'onoma
    Metropolitana, Iztapalapa, Apartado 55-534, M\'exico D.F. 09340,
    M\'exico.}

\begin{abstract}
  We present an algorithmic proof of the \CD theorem on generalized
  real scalar product spaces with arbitrary signature.  We use
  Clifford algebras to compute the factorization of a given \ort
  transformation as a product of reflections with respect to
  hyperplanes. The relationship with the Cartan-Dieudonn\'e-Scherk
  theorem is also discussed in relation to the minimum number of
  reflections required to decompose a given orthogonal transformation.
\end{abstract}
\end{frontmatter}

\section{Introduction}

The Cartan-Dieudonn\'e theorem is a fundamental result in the theory
of metric vector spaces.  It states that every orthogonal
transformation is the composition of reflections with respect to
hyperplanes. The classical proofs of the Cartan-Dieudonn\'e theorem
use induction on the dimension of the vector space and are not
constructive. See \cite{AR} and \cite{ST}.  Recently, Uhlig \cite{U1}
presented a constructive proof of the Cartan-Dieudonn\'e theorem for
the case of vector spaces with a positive definite inner product, and
also a constructive proof of a weaker version of the theorem for
generalized scalar product spaces of signature $(p,q)$
\cite[Thm.3]{U1}.

The matrix representation of a reflection with respect to a hyperplane
is called a Householder matrix \cite{U1}. The analogues of Householder
transformations on spaces with a non-degenerate bilinear or
sesquilinear form are studied in \cite{MMT}. Householder matrices are
also used in Gallier's book \cite[Ch. 7]{Gal}, which discusses the
Cartan-Dieudonn\'e theorem for linear and affine isometries, including
applications to QR decomposition.

In the present paper we present an alternative proof of the
Cartan-Dieudonn\'e theorem for generalized real scalar product spaces
of arbitrary signature. The proof yields an algorithm for the
factorization of a given orthogonal transformation as a product of
reflections with respect to hyperplanes. This work is a generalization
to spaces of arbitrary signature of a previous one \cite{AAR}, where
we provide an algorithmic proof of the Cartan-Dieudonn\'e theorem in
$\mathbb{R}^n$ valid over the fields $\mathbb{Q}$, $\mathbb{R}$ or
$\mathbb{C}$.

In the theory of Clifford algebras there is an alternative way to find
the image of a vector under a reflection with respect to a
hyperplane. This is done using vector multiplications under the rules
of the Clifford algebras.  This is the method that we propose for the
computations involved in the factorization algorithm. Aditionally, our
approach produces an alternative way for calculating the Householder
matrices with respect to orthogonal bases, but we don't use these
matrices in the  numerical examples.

In section 2 we present some results about vector spaces with
non-degenerate bilinear forms. We include some results about Artinian
spaces because they are important for the development of our proof of
the Cartan-Dieudonn\'e theorem. In section 3 we present the proof for
the case of spaces with a symmetric non-degenerate bilinear form of
arbitrary signature. In section 4 we propose the use of Clifford
algebras as a computational tool to obtain the reflections that give
the factorization of a given orthogonal transformation. In section 5
we present some examples of the factorization of orthogonal
matrices. Finally, in section 6 some conclusions are given, including
a comment about the relationship with the Cartan-Dieudonn\'e-Scherk
theorem, and the minimum number of reflections required to decompose a
given orthogonal transformation.

\section{Generalized scalar product spaces}

In this section we present some basic results concerning real vector
spaces equipped with a non-degenerate symmetric bilinear form. We call
such spaces (generalized) scalar product spaces. They are also known
as metric vector spaces or real orthogonal spaces. The proofs of the
results presented here can be found in \cite{ST}.

\begin{defn}
  Let $\mathcal{X}$ be a real vector space and let
  $\mathcal{B}:\mathcal{X} \times \mathcal{X}\rightarrow \mathbb{R}$
  be a map that satisfies the conditions

  (B1) (Bilinearity) For all $v,v^{\prime },w,w^{\prime }\in
  \mathcal{X}$ and for all $\lambda \in \mathbb{R}$
\begin{equation*}
  \mathcal{B} \left( \lambda v + v^{\prime }, w\right) = \lambda
  \mathcal{B} \left( v, w \right) + \mathcal{B} \left( v^{\prime }, w \right),
\end{equation*}
and
\begin{equation*}
  \mathcal{B}( v, \lambda w + w^{\prime }) =\lambda \mathcal{B}\left(
    v,w\right) +\mathcal{B}\left( v, w^{\prime }\right).
\end{equation*}

(B2) (Symmetry) For all $v,w\in \mathcal{X}$, \quad $\mathcal{B}\left(
  v, w \right) = \mathcal{B} \left( w, v \right)$.

(B3) (Non-degeneracy) For each non-zero $v$ in $\mathcal{X}$ there exists $w$
in $\mathcal{X}$ such that $\mathcal{B}\left( v,w\right) \neq 0$.

Then the pair $\left( \mathcal{X}, \mathcal{B}\right) $ is said to be
a (generalized) \textsl{scalar product space}.
\end{defn}

\begin{defn}
  Let $\mathcal{X}$ be a real vector space and $\mathcal{B}$ a
  bilinear form on $\mathcal{X}$.

\begin{enumerate}
\item The vectors $u,v$ in $\mathcal{X}$ are orthogonal if
  $\mathcal{B} \left( u,v \right) =0$.

\item A vector $u$ in $\mathcal{X}$ is called isotropic if
  $\mathcal{B} \left( u, u \right) =0$.

\item We say that $u$ is invertible if $u$ is not isotropic, that is,
  if $\mathcal{B} \left( u, u \right) \neq 0$. (This terminology will
  be justified in section 4.)

\item Let $\mathcal{W}$ and $\mathcal{V}$ be vector subspaces of
  $\mathcal{X}$.  We say that $\mathcal{W}$ and $\mathcal{V}$ are
  orthogonal if $\mathcal{B} \left( v, u \right) =0$, for all $u \in
  \mathcal{V}$ and $v\in \mathcal{W}$.

\item A subspace $\mathcal{V}$ of $\mathcal{X}$ is called null
  subspace if $\mathcal{B} \left( v, u \right) =0$, for all $u, v \in
  \mathcal{V}$.

\item Let $\mathcal{W}$ be a subspace of $\mathcal{X}$. The orthogonal
  complement of $\mathcal{W}$ is the subspace $\mathcal{W}^{\perp } =
  \left\{ u \in \mathcal{X} \mid \mathcal{B}( v,u) =0, \mathrm{\ for\
      all\ } v \in \mathcal{W} \right\}$.

\item Let $\mathcal{W}$ be a subspace of $\mathcal{X}$. We say that
  $\mathcal{W}$ is a non-degenerate subspace, relative to
  $\mathcal{B}$, if the restriction of $\mathcal{B}$ to $\mathcal{W}
  \times \mathcal{W}$ is non degenerate.
\end{enumerate}
\end{defn}

Note that the vector subspace generated by an isotropic vector  $u$ is a
null subspace of $\mathcal{X}$.

The next proposition states when a subspace and its orthogonal complement
decompose the space $\mathcal{X}$ as a direct sum.

\begin{prop}
  Let $\left( \mathcal{X},\mathcal{B}\right) $ be a generalized scalar
  product space of dimension $n$ and let $\mathcal{W}$ be a subspace
  of $\mathcal{X}$.  Then $\mathcal{X}=\mathcal{W}\oplus
  \mathcal{W}^{\perp }$ if and only if $\mathcal{W} $ is
  non-degenerate. That is, the space $\left( \mathcal{W}, \mathcal{B}
    \mid _{\mathcal{W} }\right) $ is non-degenerate, where
  $\mathcal{B}\mid _{\mathcal{W}}$ is the restriction of $\mathcal{B}$
  to the subspace $\mathcal{W}$. In particular, if $a\in \mathcal{X}$
  we have $\mathcal{X} = \mathbb{R}a\mathbb{\oplus } \left( \mathbb{R}
    a \right)^{\perp }$ if and only if $\mathcal{B}( a,a)
  \neq 0$, where $\mathbb{R}a$ denotes the subspace generated by $a$.
\end{prop}

For a proof see \cite[Prop. 149.1]{ST}.

In the remaider of this article $\left( \mathcal{X},\mathcal{B}\right)
$ will denote a (generalized) scalar product space of dimension $n$
over the field of real numbers.

Let $e=\left\{ e_1,e_2,\ldots,e_n\right\} $ be an ordered basis of the
vector space $\mathcal{X}$. For each pair of indices $i,j$ in $\{1, 2,
\ldots,n\}$ let $a_{i,j} := \mathcal{B}\left( e_i, e_j\right).$ The
matrix $A=[ a_{i,j}] $ is called the matrix of $\mathcal{B}$ relative
to the basis $e$. It describes the bilinear form $\mathcal{B }$ in the
following way. Let $v,w$ in $\mathcal{X}$ and let $x = \left[ v\right]
_e,\ y = \left[ w \right] _e$ be the coordinate (column) vectors of $v$
and $w$ with respect to the basis $e$. Then $\mathcal{B} \left( v,
  w\right) =x^tAy.$ Since $\mathcal{B}$ is non-degenerate and
symmetric, we see that $A$ is a non-singular symmetric matrix.

\begin{prop}
  Let $\left( \mathcal{X},\mathcal{B}\right)$ be a generalized scalar product
  space. Then there exist an ordered basis $e^{\ast } = \left(
    e_{1}^{\ast }, \ldots, e_{n}^{\ast } \right) $ of $\mathcal{X}$
  and nonnegative integers $p$ and $q$, with $p+q=n$, such that

\begin{enumerate}
\item $\mathcal{B} \left( e_{j}^{\ast }, e_{j}^{\ast } \right) =1$ for
  $j = 1, 2, \ldots, p.$

\item $\mathcal{B} \left( e_{j}^{\ast }, e_{j}^{\ast } \right) =-1$
  for $j = p+1, p+2, \ldots, p+q.$

\item $\mathcal{B} \left( e_{i}^{\ast }, e_{j}^{\ast }\right) =0$ for
  $i \neq j$.

\smallskip
  This means that the basis $e^{\ast } $ diagonalizes the matrix
  associated with the bilinear form $\mathcal{B}$.

\item The number of elements that satisfy $\mathcal{B} \left(
    e_{j}^{\ast }, e_{j}^{\ast } \right) =1$ is independent of the
  basis that diagonalizes the bilinear form $\mathcal{B}$.
\end{enumerate}
\end{prop}

For a proof see \cite[Prop. 159.1]{ST}.

The basis $e^{\ast }$ is called an orthonormal basis with respect to
$\mathcal{B}$. If $q=0$ then the space $\left( \mathcal{X},
  \mathcal{B}\right) $ is called positive definite.

Let $p$ and $q$ be nonnegative integers such that $n=p+q$. The
bilinear form $\mathcal{B}^{\ast }$ on the space $\mathbb{R}^{n}$
defined by
\begin{equation*}
  \mathcal{B}^{\ast } \left( x,y \right)
  = \sum_{i=1}^{p} x_{i} y_{i} - \sum_{i=p+1}^{p+q} x_{i} y_{i},
\end{equation*}
where $x=\left( x_{1}, x_2, \ldots , x_{p+q} \right) $ and $y = \left(
  y_{1}, y_2, \ldots , y_{p+q} \right) $, is symmetric and
non-degenerate.  The generalized scalar product space $\left(
  \mathbb{R}^{n}, \mathcal{B}^{\ast } \right) $ is denoted by
$\mathbb{R}^{p,q}$.

Any generalized scalar product space $\left( \mathcal{X},
  \mathcal{B}\right) $ that satisfies conditions 1,2 and 3 of the
previous proposition is isomorphic to the space $\mathbb{R}^{p,q}$.
Any such space $\left( \mathcal{X} , \mathcal{B} \right) $ is said to
have characteristic or signature $\left( p, q\right) $, and it is
usually identified with $\mathbb{R}^{p,q}$.  See
\cite[Thm. 177.1]{ST}.


\subsection{The orthogonal group}

Among the linear operators on $\left( \mathcal{X},\mathcal{B}\right) $
the most interesting are clearly those that preserve the bilinear
form.

\begin{defn}
  Let $T:\mathcal{X}\rightarrow \mathcal{X} $ be a linear
  operator. Then $T$ is an orthogonal transformation if and only if
\begin{equation*}
  \mathcal{B}\left( Tv,Tw\right) =\mathcal{B}\left( v,w\right), \qquad v,w \in
  \mathcal{X}.
\end{equation*}
The set of all the orthogonal transformations is a group, called the
orthogonal group of $\left( \mathcal{X}, \mathcal{B}\right)$, and
denoted by $\mathcal{O}\left( \mathcal{X}\right)$.
\end{defn}

The group $\mathcal{O}\left( \mathbb{R} ^{p,q}\right) $ can be
considered as the set of invertible $n\times n$ matrices $Q$ that
satisfy $Q^t A Q = A$, where $A$ is the matrix associated with the
bilinear form $\mathcal{B}$ with respect to the canonical basis of
$\mathbb{R}^{p,q}$. We denote $\mathcal{O} \left( \mathbb{R} ^{p,q}
\right)$ by $\mathcal{O} ( p, q )$.

\begin{defn}
  $\mathcal{SO} \left( p,q\right) := \left\{ Q \in \mathcal{O} \left(
      p, q\right) \mid \det \left( Q\right) =1\right\} $ is the group
  of special orthogonal transformations or rotations of
  $\mathbb{R}^{p,q} $.
\end{defn}

If $a$ is an invertible vector ( $\mathcal{B}\left( a,a\right) \neq
0$) then the subspace $\left( \mathbb{R}a\right) ^{\perp }$ has
dimension $n-1$ and it is called the hyperplane associated with
$a$. In this case, every $v$ in $\mathcal{X}$ has a unique
representation of the form $v = \lambda a+b$, with $b \in \left(
  \mathbb{R}a\right) ^{\perp }$ and $\lambda \in \mathbb{R}$. The next
proposition shows that the linear transformation $\varphi _{a} :
\mathcal{X} \rightarrow \mathcal{X}$, defined by $\varphi _{a} \left(
  v \right) = -\lambda a + b$ is orthogonal. It is called the
reflection with respect to the hyperplane $\left( \mathbb{R}a\right)
^{\perp}$. For the sake of convenience we denote $\left(
  \mathbb{R}a\right) ^{\perp }$ by $H_{a}$.

\begin{prop}
  Let ${\mathcal{W}}$ be a non-degenerate subspace of
  ${\mathcal{X}}$. Define the linear map $T : \mathcal{X} =
  \mathcal{W} \oplus \mathcal{W} ^{\perp } \rightarrow \mathcal{X}$ by
  $T \left( v \right) = x-y$, where $v = x+y,$, with $x$ in
  $\mathcal{W}$ and $y$ in $\mathcal{W}^{\perp }$. Then $T$ is an
  orthogonal transformation.
\end{prop}

\begin{lem}
  Let $a,b$ be invertible vectors such that $\mathcal{B}\left( a, a
  \right) = \mathcal{B} \left( b, b\right)$. Then there exists a
  linear map $\varphi$ such that $\varphi \left( a \right) = b$ and
  $\varphi $ is either the reflection with respect to a hyperplane or
  the composition of two reflections with respect to hyperplanes.
\end{lem}

\textit{Proof:} If $\mathcal{B} \left( a, a \right) =\mathcal{B}
\left( b, b\right)$ and $a, b$ are invertible then $a+b$ and $a-b$ are
orthogonal.  This follows from $\mathcal{B} \left( a+b, a-b \right) =
\mathcal{B}\left( a,a\right) -\mathcal{B} \left( b, b\right).$

We have $\mathcal{B}\left( a+b, a+b \right) +\mathcal{B }\left( a-b,
  a-b\right) = 4 \mathcal{B} \left( a, a \right) \neq 0.$ We deduce
that either $a+b$ is invertible, or $a-b$ is invertible, because at
least one of the summands in the above equation must be non-zero.

If $a-b$ is invertible then $\varphi _{a-b} : \mathcal{X} \rightarrow
\mathcal{X}$ is a reflection and
\[
\varphi _{a-b} \left( a \right) = \varphi _{a-b} \left( \frac 12\left(
    a - b \right) + \frac 12 \left( a + b\right) \right) = -\frac 12
\left( a - b\right) + \frac 12 \left( a + b \right) =b.
\]

If $a+b$ is invertible, then $\varphi _{a+b} : \mathcal{X}\rightarrow
\mathcal{X} $ is a reflection and
\begin{eqnarray*}
\varphi _b\varphi _{a+b}\left( a\right)& =& \varphi _b\varphi _{a+b}\left(
\frac { a-b}{2} +\frac { a+b}{2} \right) \\
& = & \varphi _b\left( \frac { a-b}{2} -\frac { a+b}{2} \right) \\
&=& \varphi _b\left( -b\right) =b.
\end{eqnarray*}
\mbox{$\blacksquare$}

The Cartan-Dieudonn\'e theorem states that every orthogonal
transformation $T$ on an $n$-dimensional generalized scalar product
space is the composition of at most $n$ reflections with respect to
hyperplanes. In the following section we will present our proof.

The main difficulty to obtain the proof of the Cartan-Dieudonn\'e
theorem appears in the case when $T(x)-x$ is a nonzero isotropic
vector for every nonisotropic vector $x$. This case leads us to
consider Artinian spaces. We present next some basic properties of
Artinian spaces that we will use in the next section.

\subsection{Artinian spaces}

The simplest example of an Artinian space is the Lorentz plane
$\mathbb{R }^{ 1,1 }$. In this plane the subspace generated by $u=(
1,1) $ is a null space of dimension $1.$

An Artinian space is a generalized scalar product space of the form
$\mathbb{R}^{ p,p }$ for some positive integer $p$. Every Artinian
space ${\mathcal{X}}$ has the following properties.

\begin{enumerate}
\item $\textrm{dim} \left( {\mathcal{X}} \right) = 2p$ is even.

\item $\mathcal{X}$ contains a null subspace of dimension $p$.

\item If $U$ is a maximal null subspace  of $\mathcal{X}$ then

\begin{enumerate}
\item $\textrm{dim} \left( U\right) =p$.

\item If $T$ is an element of $\mathcal{O} \left( {\mathcal{X}}
  \right) $ such that $T \left( U \right) =U $ then $T$ is a rotation,
  that is, $\det(T)=1$.

\end{enumerate}
\end{enumerate}

For our purposes, the main property of Artinian spaces is the
following lemma \cite[Prop. 247.1, Lemma 249.2]{ST}.

\begin{lem}
  Let $T$ be an element of $\mathcal{O}(p,q)$ such that $T(x) -x$ is a
  nonzero isotropic vector for every nonisotropic vector $x$. Then

\begin{enumerate}
\item $p=q$ and $2p$ is a multiple of 4.

\item $T$ is a rotation with fixed space $U$, where $U$ is a maximal
  null subspace and \quad Im$\left( T-I\right) =U=Ker\left(
    T-I\right). $
\end{enumerate}
\end{lem}

This result, combined with induction over the dimension of the space,
is used in the proof of the Cartan-Dieudonn\'e theorem presented in
\cite[Thm. 254.1]{ST}. In the next section we present an alternative
proof based on an algorithm to decompose a given orthogonal
transformation as product of reflections. In section 4.2 we present
 an explict formula for calculating the matrix representations of
the reflections.

\section{An alternative proof of the Cartan-Dieudonn\'e theorem}

In order to simplify the notation, from here on we write $uv =
{\mathcal{B}} \left( u , v \right)$, whenever $\mathcal{B}$ is the
symmetric non-degenerate bilinear form on the space
$\mathbb{R}^{p,q}$, and $u$ and $v$ are in $\mathbb{R}^{p,q}$. The
subspace generated by the vectors $u_1, u_2, \ldots, u_k$ is denoted
by $\langle u_1, u_2, \ldots, u_k \rangle$. Recall that, if $a$ is a
vector such that $a^2 \neq 0$ then $\varphi_a$ denotes the reflection
with respect to the hyperplane $H_a$ of all vectors that are
orthogonal to $a$.

The next lemma may be considered as a weak version of the
Cartan-Dieudonn\'e theorem. The proof is very similar to Uhlig's proof
of the analogous result \cite[Thm. 3]{U1}.

\begin{lem}
  Every orthogonal transformation on the space $\mathbb{R}^{p,q}$ can
  be expressed as the composition of at most $2n$ reflections with
  respect to hyperplanes, where $n=p+q$.
\end{lem}

\textit{Proof:} Let $T \in {\mathcal{O}} \left( p, q \right)$.  Let
$\left\{ w_{1}, w_{2}, \ldots , w_{n} \right\} $ be an orthogonal
basis for $\mathbb{R}^{p,q}$ such that $w_{i}^{2} \neq 0$ for $i = 1,
2, \ldots ,n.$ Define  $V_j = \langle w_j, w_{j+1},
\ldots, w_n \rangle$ for $j = 1, 2, \ldots, n$.

Consider the vectors $T \left( w_{1} \right) $ and $w_{1}$. Define the
linear function $\varphi _{1} : \mathbb{R}^{p,q} \rightarrow
\mathbb{R}^{p,q}$ by
\begin{equation*}
\varphi _{1} =%
\begin{cases}
  I, & \mathrm{\ if\ } T \left( w_{1} \right) = w_{1}, \\
  \varphi _{c_{1}}, & \mathrm{\ if\ } T \left( w_{1} \right) \neq
  w_{1} \mathrm{\ and \ } \left( T \left( w_{1} \right) -w_{1}
  \right)^{2} \neq 0, \\
  \varphi _{w_{1}} \varphi _{d_{1}}, & \mathrm{\ if\ } T \left(
    w_{1} \right) \neq w_{1} \mathrm{\ and\ } \left( T \left(
      w_{1} \right) - w_{1} \right) ^{2}=0,
\end{cases}
\end{equation*}
where $c_{1} = T \left( w_{1} \right) - w_{1}$ and $d_{1} = T \left(
  w_{1} \right) + w_{1}$. By Lemma 2.8 it is easy to see that $\varphi
_{1} T \left( w_{1} \right) = w_{1}$, and $\varphi _{1} T \left( V_2
\right) \subseteq V_2$.

Consider now the vectors $\varphi _{1} T \left( w_{2} \right)$ and
$w_{2}$.  Define the linear function $\varphi _{2} : \mathbb{R}^{p,q}
\rightarrow \mathbb{R}^{p,q}$ by
\begin{equation*}
\varphi _{2} =%
\begin{cases}
  I, & \text{\ if\ } \varphi _{1} T \left( w_{2} \right) = w_{2}, \\
  \varphi _{c_{2}}, & \text{\ if\ } \varphi _{1} T \left( w_{2}
  \right) \neq w_{2} \text{\ and\ } \left( \varphi _{1} T \left( w_{2}
    \right) -w_{2} \right) ^{2} \neq 0, \\ \varphi _{w_{2}} \varphi
  _{d_{2}}, & \text{\ if\ } \varphi _{1} T \left( w_{2} \right) \neq
  w_{2} \text{\ and\ } \left( \varphi _{1} T \left( w_{2} \right) -
    w_{2} \right) ^{2} = 0,
\end{cases}
\end{equation*}
where $c_{2} = \varphi _{1} T \left( w_{2} \right) -w_{2}$ and $d_{2}
= \varphi_{1} T \left( w_{2} \right) +w_{2}$. We know that $\varphi
_{2} \varphi_{1} T \left( w_{2} \right) = w_{2}$, and since $c_{2},
d_{2} \in V_2$, we see that $\varphi _{2} \varphi _{1} T
\left( w_{1} \right) = w_{1}$. Therefore we have $\varphi _{2} \varphi
_{1} T \left( w_{i} \right) = w_{i}$ for $i = 1,2$, and $\varphi _{2}
\varphi _{1} T \left( V_3 \right) \subseteq V_3$.

Consider now $\varphi _{2} \varphi _{1} T \left( w_{3} \right) $ and
$w_{3}$.  Define the linear function $\varphi _{3} : \mathbb{R}^{p,q}
\rightarrow \mathbb{R}^{p,q}$ by
\begin{equation*}
\varphi _{3} =
\begin{cases}
  I, & \text{\ if\ } \varphi _{2} \varphi _{1} T \left( w_{3} \right)
  = w_{3}, \\
  \varphi _{c_{3}}, & \text{\ if\ } \varphi _{2} \varphi _{1} T \left(
    w_{3} \right) \neq w_{3} \text{\ and\ } \left( \varphi _{2}
    \varphi_{1} T \left( w_{3} \right) - w_{3} \right) ^{2} \neq 0, \\
  \varphi _{w_{3}}\varphi _{d_{3}}, & \text{\ if\ }\varphi_2 \varphi
  _{1} T\left( w_{3}\right) \neq w_{3}\text{\ and\ }\left( \varphi
    _{2}\varphi _{1}T\left( w_{3}\right) -w_{3}\right) ^{2}=0,
\end{cases}
\end{equation*}
where $c_{3} = \varphi _{2} \varphi _{1} T \left( w_{3} \right)
-w_{3}$ and $d_{3} = \varphi _{2} \varphi _{1} T \left( w_{3} \right)
+ w_{3}$. We know that $\varphi _{3} \varphi _{2} \varphi _{1} T
\left( w_{3} \right) = w_{3}$, and since $c_{3}, d_{3} \in V_3$, 
we can show that $\varphi _{3} \varphi_{2} \varphi _{1} T \left(
  w_{i} \right) = w_{i}$ for $i = 1, 2, 3$. Therefore we have
$\varphi_{3} \varphi_{2} \varphi_{1} T \left( w_{i} \right) = w_{i}$
for $i = 1, 2, 3$, and $\varphi _{3} \varphi _{2} \varphi _{1} T
\left( V_4 \right) \subseteq V_4$.

We introduce the notation $\Phi_k = \varphi _{k} \varphi _{k-1} \cdots
\varphi_{1}$ for $k \ge 1$.  Then, following the procedure used above,
we can get orthogonal transformations $\varphi _{1}, \varphi _{2},
\ldots , \varphi _{n}$ such that $\Phi _{k} T \left( w_{i} \right) =
w_{i}$ for $i = 1, 2,\ldots ,k$, and $\Phi _{k} T \left( V_{k+1}
\right) \subseteq V_{k+1}$.

The orthogonal transformations $\varphi_j$ are defined by
\begin{equation*}
\varphi _{k+1}=
\begin{cases}
  I, & \text{\ if\ }\Phi _{k} T\left( w_{k+1}\right) =w_{k+1}. \\
  \varphi _{c_{k+1}}, & \text{\ if\ } \Phi _{k} T \left( w_{k+1}
  \right) \neq w_{k+1} \text{\ and\ } \left( \Phi _{k} T \left(
      w_{k+1} \right) -w_{k+1} \right)^{2} \neq 0, \\
  \varphi _{w_{k+1}}\varphi _{d_{k+1}}, & \text{\ if\ }\Phi _{k}
  T\left( w_{k+1}\right) \neq w_{k+1} \text{\ and \ }\left( \Phi
    _{k}T\left( w_{k+1}\right) -w_{k+1}\right) ^{2}=0,
\end{cases}
\end{equation*}
where $c_{k+1} = \Phi _{k}T\left( w_{k+1}\right) -w_{k+1},$ and
$d_{k+1} = \Phi _{k} T \left( w_{k+1} \right) +w_{k+1}$.

Therefore we have $\Phi _{n}T\left( w_{i}\right) =w_{i}$ for
$i=1,2,\ldots ,n$, and thus $\Phi _{n}T=I$, and then
\begin{equation*}
  T =\Phi_n^{-1}=\varphi _{1}^{-1} \varphi _{2}^{-1} \cdots \varphi _{n}^{-1}.
\end{equation*}
Since each $\varphi _{k}^{-1} $ is either the identity, a reflection
with respect to a hyperplane, or the composition of two reflections
with respect to hyperplanes, we see that $T$ is the composition of at
most $2n$ reflections with respect to hyperplanes and this completes
the proof. \mbox{$\blacksquare$}

Reviewing the ideas used in the proof of the previous lemma, we see
that one way to reduce the number of reflections needed to factor $T$
is the following. Suppose we have found $\varphi _{1}, \varphi _{2},
\ldots ,\varphi _{k}$, reflections with respect to hyperplanes, such
that $\Phi _{k} T \left( w_{i} \right) = w_{i}$ for $i = 1, 2, \ldots
,\ell$, and $\Phi _{k} T \left(V_{\ell+1} \right) \subseteq V_{\ell
  +1}$. If there exists $ w_{j} \in \{w _{\ell+1}, w_{\ell+2}, \ldots,
w_n\}$ such that $\Phi _{k} T \left( w_{j} \right) = w_{j}$, or $\Phi
_{k} T \left( w_{j} \right) \neq w_{j}$ and $\left( \Phi _{k} T \left(
    w_{j} \right) - w_{j} \right) ^{2} \neq 0,$ then we can reorder
the elements to force $j= \ell +1$, and then, using the construction
of the $\varphi_i$ of the previous lemma we would have that $\varphi
_{k+1}$ is either the identity or a reflection with respect to a
hyperplane.

Lemma 2.9 tells us under what conditions it can happen that no element
of $\left\{ w_{\ell+1}, w_{\ell+2}, \ldots , w_{n} \right\} $
satisfies the conditions described above, and in such situation we can
not assure that $\varphi _{k+1}$ is the identity or a reflection with
respect to a hyperplane.

In order to get a proof of the Cartan-Dieudonn\'e theorem, we must
find an algorithm for the construction of the orthogonal
transformations $\varphi_j$ that avoids in some way reaching a
situation where the hypothesis of Lemma 2.9 is satisfied. This can be
done by introducing an additional reflection in the way we describe
next.

Let $T \in {\mathcal{O}}(p,q)$ and suppose that $T(x)-x$ is a nonzero
isotropic vector for every nonisotropic $x$. Then, by Lemma 2.9 we
must have $p=q$, $n=p+q$ is a multiple of 4, and $T$ is a rotation,
that is, det$(T)=1$.

Let $y$ be an invertible element of $\mathbb{R}^{p,q}$ and let
$\varphi_{y}$ be the reflection with respect to the hyperplane
$H_{y}$. Define $S=\varphi_{y}T$. Since det$(S)=-1$ we see that $S\neq
I$ and $S$ is not a rotation. Therefore $S$ does not satisfy the
hypothesis of Lemma 2.9.

Let $\left\{ w_{1}, w_{2}, \ldots ,w_{n}\right\} $ be an orthogonal
basis for $\mathbb{R}^{p,q}$. We can reorder the elements of the basis
so that either $S \left( w_1 \right) = w_1$ or $S \left( w_{1} \right)
-w_{1}$ is an invertible vector. Then we can find $\varphi _{1}$ that
is either the identity or a reflection with respect to a hyperplane,
and satisfies $\varphi _{1} S \left( w_{1} \right) = w_{1}$ and
$\varphi _{1} S \left( V_2 \right) \subseteq V_2$.

Since the dimension of $V_2$ is not a multiple of 4, the orthogonal
transformation $\varphi_1 S$ restricted to $V_2$ does not satisfy the
hypothesis of Lemma 2.9. Therefore there exists $j$ such that $2 \le j
\le n $ and either $\varphi _{1} S \left( w_{j} \right) =w_{j}$ or
$\varphi _{1} S \left( w_{j}\right) -w_{j}$ is an invertible
vector. Reordering the basis of $V_2$ if necessary, we can suppose
that $j=2$. Then we can find $\varphi _{2}$ that is either the
identity or a reflection with respect to a hyperplane and satisfies
$\varphi _{2}\varphi _{1} S \left( w_{i} \right) =w_{i}$ for $i=1,2$
and $\varphi _{2} \varphi _{1} S \left( V_3 \right) \subseteq V_3$.

Proceeding in the same way we can get $\varphi _{3}$ and $\varphi
_{4}$, that are either the identity or reflections, that satisfy
$\varphi _{4} \varphi_{3} \varphi _{2} \varphi _{1} S \left( w_{i}
\right) =w_{i}$ for $i=1, 2, 3, 4$, and $\varphi _{4} \varphi _{3}
\varphi _{2} \varphi _{1} S \left( V_5 \right) \subseteq V_5$.

Consider now the composition $\Phi_4 S$, where $\Phi_4 = \varphi_{4}
\varphi_{3} \varphi_{2} \varphi_{1}$.  (Recall that we defined $\Phi_k
= \varphi_{k} \varphi_{k-1} \cdots \varphi_{2} \varphi _{1}$ for $k
\ge 1$). There are two possible cases:

\begin{enumerate}
\item $\varphi _{i}\neq I$ for $i=1,2,3,4$.

  In this case we have $\det ( \Phi_4 S) =-1$ and hence $\Phi_4 S$
  restricted to $V_5$ does not satisfy the hypothesis of Lemma
  2.9. Therefore we can find orthogonal transformations $\varphi_j$,
  for $j=5,6,7,8$ such that $\Phi_8 S \left( w_{i} \right) =w_{i}$ for
  $i=1, 2, \ldots ,8$ and $\Phi_8 S \left(V_9 \right) \subseteq
  V_9$. Then, for the composition $\Phi_8 S$ we have again the same
  two possible cases that we had for $\Phi_4 S $, but now considering
  the maps $\varphi _{i}$ for $i = 5, 6, 7, 8$.

\item $\varphi _{i}=I$ for at least one $i$, with $1 \le i \le 4$.

  In this case $\dim \left(V_5 \right)$ is a multiple of 4 and it is
  possible that $\det \left( \Phi_4 S\right) =1$.

  If $\Phi_4 S$ restricted to the space $V_5$ does not satisfy the
  hypothesis of Lemma 2.9 then we can find $\varphi _{5}$, which is
  either the identity or a reflection such that $\Phi_5 S\left(
    w_{i}\right) =w_{i}$ for $i = 1, 2, \ldots , 5$ and $\Phi_5 S
  \left( V_6 \right) \subseteq V_6$. Notice that, since at least one
  of the $\varphi _{i}$ is the identity, for $1 \le i \le 4 $, the
  number of reflections in the composition $\Phi_5$ is at most equal
  to 4.

  If $\Phi_4 S$ restricted to the space $V_5$ satisfies the hypothesis
  of Lemma 2.9 then we take an invertible vector $z$ in $V_5$ and form
  the composition $\varphi_z \Phi_4 S$. This map restricted to $V_5$
  can not satisfy the hypothesis of Lemma 2.9, and consequently, we
  can find $\varphi_{5}$, which is either the identity or a reflection
  and satisfies $\varphi _{5}\varphi_z \Phi_4 S \left( w_{i} \right)
  =w_{i}$ for $i= 1, 2, \ldots ,5$ and $\varphi _{5} \varphi_z \Phi_4
  S \left( V_6 \right) \subseteq V_6$. Since at least one of the
  $\varphi _{i}$ is the identity, for $1 \le i \le 4 $, the number of
  reflections in the composition $\varphi _{5} \varphi_z \Phi_4$ is at
  most equal to 5.
\end{enumerate}

Applying the procedure described above we see that for each $\ell$
such that $1 \le \ell \le n$, we can find reflections with respect to
hyperplanes $\varphi_1, \varphi_2, \ldots, \varphi_s$ such that
$\Phi_s S \left( v_i \right) = v_i$ for $1 \le i \le \ell$, where
$\left \{ v_{1}, \ldots , v_{n} \right\} $ is a reordering of the
orthogonal basis $\left\{ w_{1}, \ldots , w_{n} \right\} $, and $s\leq
\ell$. This last inequality is very important.

In particular, for $\ell=n$ we get $\Phi_s S = I$, with $s \le n$. We
claim that the case $s=n$ is not possible. If $s=n$ then $\det \left(
  \Phi_s \right)=(-1)^n=1$, because $n$ is a multiple of 4. On the
other hand, $\det \left( \Phi_s \right) \det ( S ) = \det(I)=1$. But
we know that $\det(S)=-1$. Therefore $s=n$ is not possible and we
conclude that $s < n$.

Since $S = \varphi_y T$, we have $\Phi_s \varphi_y T =I$ and therefore
$T$ is the composition of at most $n$ reflections with respect to
hyperplanes.

We have proved the following result.

\begin{lem}
  Let $T$ be an element of ${\mathcal{O}} (p, q)$. If $T\left( x
  \right) -x$ is a nonzero isotropic vector for every nonisotropic
  vector $x$ then $T$ is the composition of at most $p+q$ reflections
  with respect to hyperplanes.
\end{lem}

\begin{lem}
 \label{lema:nuevo}
  Let $T$ be an element of ${\mathcal{O}} (p, q)$. If there exists a
  basis $\left\{ w_1, \ldots, w_{p+q} \right\}$, where all the elements
    are nonisotropic, such that
\[
T\left( w_{i} \right) - w_{i}
\]
is a nonzero isotropic vector for $i = 1, \ldots, p+q$, then $T$ is
the composition of at most $p+q$ reflections with respect to
hyperplanes.
\end{lem}

\textit{Proof:} We can proceed as follows.

\begin{description}
\item[Step~1] In each step we have to deal with three possible
  cases. In particular, here we have:
\begin{enumerate}
\item There exists a nonzero and nonisotropic element $v_{1} \in
  \mathbb{R}^{p,q}$ such that $T \left( v_{1} \right) - v_{1} = 0$.

\item There exists a nonzero and nonisotropic element $v_{1}$ of the
  basis such that $T \left( v_{1} \right) - v_{1} \neq 0$ and $\left(
    T \left( v_{1} \right) - v_{1} \right) ^{2}\neq 0$ (\emph{i.e.} $T
  \left( v_{1} \right) - v_{1} $ is nonisotropic).

\item For every nonzero and nonisotropic $x \in \mathbb{R}^{p,q}$ we
  have that $T \left( x \right) - x \neq 0$ is isotropic.
\end{enumerate}

In case 3), from Lemma 3.2 we obtain that $T$ is the composition of at
most $p+q$ reflections with respect to hyperplanes.

In cases 1) or 2), we can find $\varphi _{1}$, which is either the
identity or a reflection, that satisfies
\[
\varphi _{1} T \left( v_{1} \right) = v_{1}
\]
and also $\varphi _{1} T \left( W_1 \right) = W_1$, where $W_1 =
\langle w_1 \rangle ^\perp$.

\item[Step~2] Now, for the orthogonal transformation $\varphi _{1} T$
  restricted to the space $W_1$, we have
\begin{enumerate}
\item There exists a nonzero and nonisotropic element $v_{2} \in W_1$
  such that $\varphi_1 T \left( v_{2} \right) - v_{2} = 0$.

\item There exists a nonzero and nonisotropic element $v_{2} \in W_1$
  such that $\varphi_1 T \left( v_{2} \right) - v_{2} \neq 0$ is
  nonisotropic.

\item For each nonzero and nonisotropic $x \in \mathbb{R}^{p,q}$ we
  have that $\varphi_1 T \left( x \right) - x \neq 0$ is isotropic.
\end{enumerate}

In case 3), we have that $\varphi_1 T = S$ and thus, from Lemma 3.2,
$S$ is the composition of at most $p+q-1$ reflections with respect to
hyperplanes.

In cases 1) or 2), we can find $\varphi _{2}$, which is either the
identity or a reflection, that satisfies
\[
\varphi_2 \varphi _{1} T \left( v_{i} \right) = v_{i},
\] 
where $i$ can be chosen as $1$ or $2$. Also, it is fulfilled that
$\varphi_2 \varphi _{1} T \left( W_2 \right) = W_2$, where $W_2 =
\langle w_1, w_2 \rangle ^\perp$.
\end{description}

By following these steps, we can end up with one of the two following
situations:

\begin{description}
\item[A] We can find an orthogonal set of nonisotropic elements
  $\left\{v_1,v_2, \ldots ,v_{p+q}\right\}$ and a finite sequence of
  linear transformations $\varphi_1,\varphi_2, \ldots, \varphi_{p+q}$,
  such that:
      \begin{itemize}
      \item $\varphi_i$ is either a reflection or the identity, for
        $i=1, 2, \ldots, p+q$.

      \item $\varphi_{l} \varphi_{l-1} \cdots \varphi_{2} \varphi_{1}
        T \left( v_i \right) = v_i$, for $i = 1, 2, \ldots, l$ and $l =
        1, 2, \ldots p+q$.

      \item $\varphi_{l} \varphi_{l-1} \cdots \varphi_{2} \varphi_{1}
        T \left( W_i \right) = W_i$, for $l = 1, 2, \ldots p+q$, where
        $W_i = \left\langle v_1, v_2, \ldots , v_i \right\rangle
        ^\bot$ for $i = 1, 2, \ldots p+q$.
      \end{itemize}

    \item[B] We can find an orthogonal set of nonisotropic elements
      $\left\{v_1, v_2, \ldots , v_{k}\right\}$ and a finite sequence
      of linear transformations $\varphi_1, \varphi_2, \ldots,
      \varphi_{k}$, where $k<p+q$, such that:
      \begin{itemize}
      \item $\varphi_i$ is either a reflection or the identity, for $i
        = 1, 2, \ldots, k$.

      \item $\varphi_{l} \varphi_{l-1} \cdots \varphi_{2} \varphi_{1}
        T \left( v_i \right) = v_{i}$, for $i = 1, 2, \ldots, l$ and
        $l = 1, 2, \ldots k$.

      \item $\varphi_{l} \varphi_{l-1} \cdots \varphi_{2} \varphi_{1}
        T \left( W_i \right) = W_i$, for $l = 1, 2, \ldots k$, where
        $W_i = \left\langle v_1, v_2, \ldots , v_i \right\rangle
        ^\bot$ for $i = 1, 2, \ldots k$.

      \item Consider the orthogonal transformation $\varphi_{k}
        \varphi_{l-1} \cdots \varphi_2 \varphi_{1} T$, restricted to
        $W_k$. Then we have that for every nonzero and nonisotropic $x
        \in W_k$ it is fulfilled that $\varphi_{k} \varphi_{l-1}
        \cdots \varphi_2 \varphi_{1} T \left( x \right) - x$ is
        nonzero and isotropic.
      \end{itemize}
\end{description}

In case \textbf{A}, we have that
\[
\varphi_{p+q} \varphi_{p+q-1} \cdots \varphi_{2} \varphi_{1} T = I.
\]
Thus $T$ is the composition of at most $p+q$ reflections.

In case \textbf{B}, we have that
\[
\varphi_{k} \varphi_{l-1} \cdots \varphi_2 \varphi_{1} T = S,
\]
restricted to $W_k$ (whose dimension is $p+q-k$). From Lemma 3.2 we
conclude that $S$ is the composition of at most $p+q-k$
reflections.

In summary, $T$ is the composition of at most $p+q$
reflections. \mbox{$\blacksquare$}

Now we are ready to prove the Cartan-Dieudonn\'e theorem.

\begin{thm}
  Let $T$ be an element of ${\mathcal{O}}(p,q)$.  Then $T$ is the
  composition of at most $p+q$ reflections with respect to
  hyperplanes.
\end{thm}

\textit{Proof:} Let $n=p+q$ and let $\left\{ w_{1}, w_{2}, \ldots ,
  w_{n} \right\} $ be an orthogonal basis for $\mathbb{R}^{p,q}$ with
$w_{i}^{2} \neq 0$ for $1 \le i\le n$.

If there exists an element $w_{j}$ of the basis such that $T \left(
  w_{j} \right) = w_{j} $ or $T \left( w_{j} \right) \neq w_{j}$ and
$\left( T \left( w_{j} \right) - w_{j} \right) ^{2}\neq 0$, then,
reordering the basis (keeping the notation $w_i$ for the basic
elements), we can get $j=1$. Thus we can find $\varphi _{1}$, which is
either the identity or a reflection, that satisfies $\varphi _{1} T
\left( w_{1} \right) = w_{1}$, and $\varphi _{1} T \left( V_2 \right)
\subseteq V_2$.

If there exists an element $w_{j} \in \left\{ w_{2}, w_{3}, \ldots ,
  w_{n} \right\} $ such that $\varphi _{1} T \left( w_{j} \right) =
w_{j}$ or $\varphi _{1} T \left( w_{j} \right) \neq w_{j}$ and $\left(
  \varphi _{1} T \left( w_{j} \right) - w_{j} \right) ^{2} \neq 0$,
reordering the basis $\left\{ w_{2}, w_{3}, \ldots , w_{n} \right\} $
we can assume that $j=2$, and then we can find $\varphi _{2},$ that is
either the identity or a reflection with respect to a hyperplane, that
satisfies $\varphi _{2} \varphi _{1} T \left( w_{i} \right) = w_{i}$
for $i = 1, 2$, and $\varphi _{2} \varphi _{1} T \left( V_3 \right)
\subseteq V_3$.

This process can be repeated until we reach either of two possible
cases. In the first case we can find $\varphi_i$, for $1 \le i \le n$,
such that $\Phi_n T = I$, and therefore $T$ is the composition of at
most $n$ reflections with respect to hyperplanes.

In the second case we find $\varphi_i$, for $1 \le i \le k$, where $1
\le k < n$, that satisfy $\Phi _{k} T \left( w_{i} \right) =w_{i}$ for
$1 \le i \le k$, $\Phi_k T \left( V_{k+1} \right) \subseteq V_{k+1}$,
and $\Phi_k T(w_j) - w_j$ is a nonzero isotropic vector for every
para $k<j\leq n$. Applying Lemma 3.3 to the
orthogonal transformation $\Phi_k T$, on the space $V_{k+1}$, we
obtain reflections with respect to hyperplanes $\tau_i$, for $1 \le i
\le s $, with $s \le n-k$, such that $\tau_s \tau_{s-1}\cdots \tau_1
\Phi_k T $ is the identity on ${\mathbb{R}}^{p,q}$. Since $s+k\le n$,
we conclude that, also in this case, $T$ is the composition of at most
$n$ reflections with respect to hyperplanes. \mbox{$\blacksquare$}

\section{Clifford algebras}

In this section we present some definitions and basic results about
Clifford algebras associated with a generalized scalar product space
$\left( \mathcal{X},\mathcal{B}\right) $ of signature $(p,q)$. We use
the notation of \cite[Ch. I]{DS}. We also describe how the Clifford
algebra structure can be used to deal with orthogonal transformations,
and in particular, with reflections with respect to hyperplanes.

\begin{defn}
  Let $\left( \mathcal{X},\mathcal{B}\right) $ be a generalized scalar
  product space of dimension $n$ and let $\mathcal{A}$ be a real
  associative algebra with identity 1 such that

  (C1) $\mathcal{A}$ contains copies of $\mathbb{R}$ and of
  $\mathcal{X}$ as linear subspaces.

  (C2) For all $v\in \mathcal{X}$ we have $v^{2}=\mathcal{B}\left(
    v,v\right).$

  (C3) $\mathcal{A}$ is generated as a ring by the copies of
  $\mathbb{R}$ and of $\mathcal{X}$, or equivalently, as a real
  algebra by $\left\{ 1\right\} $ and $\mathcal{X}$.

  Then $\mathcal{A}$ is called a real \textit{Clifford algebra} for
  $\left( \mathcal{X},\mathcal{B}\right) $ and it is denoted by
  $\mathcal{A}= \mathcal{C} \left( \mathcal{X} \right)$.
\end{defn}

Note that (C2) links the multiplication in the algebra with the
bilinear form on $\mathcal{X}$.


\subsection{Bases for Clifford algebras}

Here we describe a basis for $\mathcal{C}\left( \mathcal{X}\right) $
in terms of an orthonormal basis for $\mathcal{X}$. Let $\left(
  \mathcal{X}, \mathcal{B}\right) $ be a generalized scalar product
space of signature $(p,q)$ and let $e = \left\{e_1, e_2, \ldots, e_n
\right\}$ be an orthonormal basis for $\mathcal{X}$. Then, by (C2) we
have

\begin{equation*}
  e_i^2={\mathcal{B}}(e_i,e_i)=
\begin{cases}
  1, & \mathrm{for\ }i=1,2,\ldots,p \\
  -1, & \mathrm{for\ } i=p+1,p+2,\ldots,p+q,
\end{cases}
\end{equation*}
and it is easy to show that
\begin{equation*}
  e_ie_j+e_je_i=0, \qquad i\ne j,
\end{equation*}
and
\begin{equation*}
  \frac 12\left( vw+wv\right) =\mathcal{B} \left( v, w \right), \qquad
  u, v \in \mathcal{X}.
\end{equation*}
Define $N=\left\{ 1, 2, \ldots, n \right\}$. Let $\beta _1, \beta_2,
\ldots , \beta _s$ be distinct elements of $N$. Then
\begin{equation*}
  e_{\beta _1} e_{\beta_2} \cdots e_{\beta _s} = \left( -1\right) ^\sigma
  e_{\alpha _1} e_{\alpha_2} \cdots e_{\alpha _s},
\end{equation*}
where $\left( \alpha_1, \alpha_2, \ldots, \alpha_s \right)$ is the
permutation of $\left( \beta_1, \beta_2, \ldots, \beta_s \right)$ that
satisfies $\alpha_1 < \alpha_2 < \cdots < \alpha_n$, and $\sigma$ is
the number of transpositions of the permutation that sends $\left(
  \alpha_1, \alpha_2, \ldots, \alpha_s \right)$ to $\left( \beta_1,
  \beta_2, \ldots, \beta_s \right)$.

Since $e_i^2=\pm 1$, it is easy to see that, up to a change of sign,
every product of basic elements, possibly with repeated factors, can
be reduced to a product of at most $n$ factors with indices in
increasing order. This implies that every element of $\mathcal{C}
\left( \mathcal{X} \right) $ can be written in the form
\begin{equation*}
  \sum\limits_A \lambda _A \, e_{\alpha _1}e_{\alpha _2} \cdots e_{\alpha _s},
\end{equation*}
where the sum runs over all the subsets $A = \left\{ \alpha_1,
  \alpha_2, \ldots, \alpha_n \right\}$ of $N$ (with $\alpha_1 <
\alpha_2 < \cdots <\alpha_n$) and $\lambda_A$ is a real
coefficient. In order to simplify the notation we define $e_A =
e_{\alpha _1} e_{\alpha _2} \cdots e_{\alpha _s}$, and we put
$e_{\emptyset} = 1$. Therefore the collection of elements of the form
$e_A$, where $A$ is a subset of $N$, is a generating set for
$\mathcal{C} \left( \mathcal{X} \right)$ and consequently,
$\textrm{dim} \left( \mathcal{C} \left( \mathcal{X} \right) \right)
\le 2^n$.

Recall that a vector $s\in \mathcal{X}$ is called invertible if
$\mathcal{B} \left( s, s \right) \neq 0$. This condition is equivalent
to $s^2 = \mathcal{B} \left( s, s \right)\ne 0$. Therefore, the
element $\left( 1/ \mathcal{B} \left( s, s\right) \right) s$ is the
inverse of $s$ with respect to the multiplication in the algebra
$\mathcal{C} \left( \mathcal{X} \right)$, and we can write $s^{-1} =
s/s^2 = s/ \mathcal{B} \left( s,s \right)$.

\begin{thm}
  Let $\mathcal{C} \left( \mathcal{X} \right) $ be a Clifford algebra
  for an $n$-dimensional generalized scalar product space $\left(
    \mathcal{X}, \mathcal{B} \right) $ with signature $\left( p, q
  \right) $ and let 

$e = \left\{ e_{1}, \ldots , e_{p}, e_{p+1}, \ldots , e_{p+q} \right\} $
 be an orthonormal basis for 
$\left( \mathcal{X}, \mathcal{B} \right) $. Then

  $\left( i\right) $ If $n$ is even then $\dim \mathcal{C} \left(
    \mathcal{X} \right) =2^{n}$ and $\{ e_{A}: A\subseteq N \} $ is a
  basis for $\mathcal{C} \left( \mathcal{X} \right)$.

  $\left( ii\right) $ If $n$ is odd and $e_{N} \notin \mathbb{R}$,
  then $\dim \mathcal{C} \left( \mathcal{X}\right) = 2^{n}$ and
  $\left\{ e_{A}: A \subseteq N \right\} $ is a basis for $\mathcal{C}
  \left( \mathcal{X } \right)$.

  $\left( iii\right) $ If $n$ is odd and $e_{N} \in \mathbb{R}$ , then
  $e_{N} = \pm 1$ and $p-q \equiv 1 \ \ {\mathrm{\ }mod\ } \left(
    4\right)$. In this case $\dim \mathcal{C} \left( \mathcal{X}
  \right) =2^{n-1}$ and $\left\{ e_{A}: A \subseteq N, \text{ \#} A
    \text{\ even\ } \right\} $ is a basis for $\mathcal{C} \left(
    \mathcal{X} \right)$.
\end{thm}

\begin{thm}
  Let $\mathcal{C} \left( \mathcal{X} \right) $ be a Clifford algebra
  for the generalized scalar product space $\left( \mathcal{X},
    \mathcal{B} \right)$ and let $W$ be a non-degenerate subspace of
  $\mathcal{X}$. Then the subalgebra of $\mathcal{C} \left(
    \mathcal{X} \right)$ generated by $W$ is a Clifford algebra for
  $W$.
\end{thm}

The proof of these results can be found in \cite{DS}, \cite{IP}, and
\cite{MR}.

A result of Chevalley states that every $n$-dimensional generalized
scalar product space $\left( \mathcal{X}, \mathcal{B}\right) $ has a
Clifford algebra of dimension $2^n$ (see \cite{DS}). If $\left(
  \mathcal{X}, \mathcal{B} \right) $ has signature $(p,q)$ then it is
isomorphic to ${\mathbb{R}}^{p,q}$. The Clifford algebra of dimension
$2^n$ of ${\mathbb{R}}^{p,q}$ is denoted by ${\mathbb{R}}_{p,q}$.

The algebra ${\mathbb{R}}_{p,q}$ is a graded vector space. Let $s$ be
an integer such that $0 \le s \le n $. Let ${\mathbb{R}}_{p,q}^s$
denote the subspace generated by $\{ e_A: A \subseteq N, \ \#A= s
\}$. It has dimension $\binom{n}{s}$. The elements of
${\mathbb{R}}_{p,q}^s$ are called $s$-vectors. Note that
\begin{equation*}
  \mathbb{R}_{p,q} = \bigoplus_{s=1}^{n} \mathbb{R}_{p,q}^{s}.
\end{equation*}
The space of $0$-vectors is generated by $\{1 \}$, and its elements
are also called scalars. The space ${\mathbb{R}}_{p,q}^1$ can be
identified with ${\mathbb{R}}^{p,q}$. The space ${\mathbb{R}}_{p,q}^n$
is generated by $\{ e_N\}$. The element $e_N$ is called pseudoscalar.

Every element $a\in \mathbb{R}_{p,q}$ can be expressed in the form
\begin{equation*}
  a =\sum\limits_{A\subseteq N}\lambda _{A}e_{A}, \qquad \lambda _{A} \in
  \mathbb{R}.
\end{equation*}
and also as
\begin{equation*}
  a =\sum_{r=0}^{n}a_{r},\qquad a_{r}\in \mathbb{R}_{p,q}^{r}.
\end{equation*}

\begin{defn}
  Let $a = \sum_{r=0}^n a_r$, where $a_r \in \mathbb{R}_{p,q}^r$. If
  $a_r = 0$ for $r > t$ and $a_t \neq 0$, then we say  that the
  grade of the multivector $a$ is $t$, and 
write $\textrm{gr} ( a ) = t.$
\end{defn}
Notice that if $\lambda \in \mathbb{R}$ then $\textrm{gr} ( \lambda a
) = \textrm{gr} ( a )$. If $s_1, \ldots, s_k$ are non-zero vectors in
$\mathbb{R}^{p,q}$, then $\textrm{gr} \left( s_1 \cdots s_k \right) =
l \leq k$ and $\textrm{dim} \langle s_1, \ldots, s_k \rangle = l$.

The multiplication in the algebra $\mathbb{R}_{p,q}$ can be used to
represent reflections with respect to a hyperplane in the space
$\mathbb{R}^{p,q}$ as follows. Let $s$ be an invertible vector in
$\mathbb{R}^{p,q}$.  Define the map $\varphi _{s}
:\mathbb{R}^{p,q}\rightarrow \mathbb{R}^{p,q}$ by
\begin{equation*}
  \varphi _{s}\left( x\right) =-sxs^{-1}, \qquad x \in \mathbb{R}^{p,q}.
\end{equation*}

In order to show that $\varphi _{s} \left( x \right)$ is a vector in
$\mathbb{R}^{p,q} $ for all $x\in \mathbb{R}^{p,q}$ we first compute
$sxs$ and get
\begin{eqnarray*}
  sxs &=&\frac{1}{2}\left( sx+xs+\left( sx-xs\right) \right) s \\
  &=&\frac{1}{2}\left( 2\mathcal{B}\left( x,s\right) +sx-xs\right) s \\
  &=&\mathcal{B}\left( x,s\right) s+\frac{1}{2}sxs-\frac{1}{2}xs^{2}.
\end{eqnarray*}
The last equation yields $sxs =2 \mathcal{B} \left( x, s \right)
s-xs^{2}$ and hence
\begin{equation}
  \label{cd01}
  -sxs^{-1} =-\frac{2 \mathcal{B} \left( x, s \right) }{s^{2}}s + x,
\end{equation}
which is  clearly a vector in $\mathbb{R} ^{p,q}$.

\begin{lem}
  Let $s$ be an invertible vector in $\mathbb{R} ^{p,q}$. Then the
  linear map $\varphi _{s}:\mathbb{R}^{p,q}\rightarrow \mathbb{R}
  ^{p,q}$ is an orthogonal transformation. Furthermore, it is the
  reflection with respect to the hyperplane $H_{s}=\left\{ x\in
    \mathbb{R}^{p,q}: \mathcal{B} \left( x,s\right) =0\right\}$.
\end{lem}

\textit{Proof:} For $x,y\in \mathbb{R}^{p,q}$ we have
\begin{eqnarray*}
  \mathcal{B}\left( \varphi _{s}\left( x\right) ,\varphi _{s}\left( y\right)
  \right) &=&\mathcal{B}\left( -sxs^{-1},-sys^{-1}\right) \\
  &=&\dfrac{1}{2}\left( \left( -sxs^{-1}\right) \left( -sys^{-1}\right)
    +\left( -sys^{-1}\right) \left( -sxs^{-1}\right) \right) \\
  &=&\frac{1}{2}s\left( xy+yx\right) s^{-1} \\
  &=&\mathcal{B}\left( x,y\right).
\end{eqnarray*}
If $x\in H_{s}$ we have $\varphi _{s} \left( x\right) = -sxs^{-1} =
-\left( 2 \mathcal{B} \left( x,s \right) /s^{2} \right)\, s+x=x$.

If $x= \lambda s$ where $\lambda$ is a real number then $\varphi _{s}
\left( x \right) = -s \left( \lambda s\right) s^{-1} = -\lambda s =
-x. $ Therefore $\varphi _{s}$ is the reflection with respect to
$H_{s}$. \mbox{$\blacksquare$}

It is easy to verify that $\varphi _{s}$ satisfies

\begin{enumerate}
\item $\varphi _{s} =\varphi _{\lambda s}$ for every nonzero real
  $\lambda$.

\item $\varphi _{s}^{-1}= \varphi _{s}$.

\item If $s_1, s_2, \ldots, s_k$ are invertible elements of
  $\mathbb{R}^{p,q}$ then
\begin{equation*}
  \varphi _{s_{1}} \varphi _{s_{2}} \cdots \varphi _{s_{k}} \left( x \right)
  = \left( -1 \right) ^{k}s_{1}s_{2} \cdots
  s_{k-1} s_{k} x s_{k}^{-1} s_{k-1}^{-1} \cdots  s_{2}^{-1} s_{1}^{-1}.
\end{equation*}
\end{enumerate}

Using the Clifford algebra multiplication the Cartan-Dieudonn\'e
theorem reads as follows.

\begin{thm}
  Let $T$ be an orthogonal transformation on $\mathbb{R}^{p,q}$. Then
  there exist invertible elements $s_1, s_2, \ldots, s_k$ in
  $\mathbb{R}^{p,q}$, with $k \le p+q$, such that
\begin{equation*}
  T\left( x\right) =\left( -1\right) ^{k} s_{1} s_{2} \cdots
  s_{k-1} s_{k} x s_{k}^{-1} s_{k-1}^{-1} \cdots  s_{2}^{-1} s_{1}^{-1}, \qquad x
  \in \mathbb{R}^{p,q}.
\end{equation*}
\end{thm}

If $T$ is an orthogonal transformation on $\mathbb{R} ^{p,q}$ and
there exist invertible multivectors $A,B$ in the Clifford algebra
$\mathbb{R}_{p,q} $ such that
\begin{equation*}
  T\left( x \right) = \pm A x A^{-1} = \pm BxB^{-1}, \qquad x \in \mathbb{R} ^{p,q},
\end{equation*}
then $A=\lambda B$ for some real $\lambda$. The proof of this result can be
found in \cite{MR}. We will illustrate this fact in the examples of the next
section.

\subsection{A matrix representation of $\varphi_s$}

Having in mind that, in many applications, matrix representations of
reflections are useful, we introduce another algebraic expression for
$\varphi_s$.

Define $\lambda =\frac{\mathcal{B}\left(x,s\right)}{x^2}$, then

\begin{equation*}
  s = \lambda x+y,
\end{equation*}

where $y$ is orthogonal to $x$, that is
\[
xy = -y x.
\]
Thus
\begin{eqnarray}
 \label{cd02}
 -s x s^{-1} & = & -\frac{1}{s^2} \left( \lambda x + y \right) x \left(
   \lambda x + y\right), \notag\\
 &=& -\frac{1}{s^2} \left( \left(\lambda^2 x^2 - y^2 \right) x + 2
   \lambda x^2 y \right), \notag\\
 &=& -\frac{1}{s^2} \left( \left( \frac{ \left( \mathcal{B} \left(
           s, x \right) \right)^2}{x^2} - y^2 \right) x + 2
   \mathcal{B} \left( s , x \right) y \right) .
\end{eqnarray}
This last result (\ref{cd02}) turns out to be useful to find a matrix
representaction of $\varphi_s$ with respect to a given orthogonal
basis $B = \left\{ w_1, w_2, \ldots, w_{p+q} \right\}$, where $w_i^2
\neq 0$ for $i = 1, 2, \ldots, p+q$. Indeed, consider $s =
\sum_{i=1}^{p+q} \alpha_{i }w_{i}$, then
\begin{eqnarray*}
  s &=& \alpha_{k }w_{k} + \sum\limits_{\substack{ i=1  \\ i\neq
      k}}^{p+q} \alpha_{i }w_{i}, \\
  &=& \frac{\mathcal{B} \left( s, w_{k} \right) }{w_{k}^{2}} w_{k} +
  \sum\limits_{\substack{ i=1  \\ i\neq
      k}}^{p+q}\frac{\mathcal{B}\left( s,w_{i}\right) }{w_{i}^{2}}
  w_{i} .
\end{eqnarray*}
Define
\[
y_{k} = \sum\limits_{\substack{ i=1 \\ i\neq k}}^{p+q} \alpha_{i}
w_{i} = \sum\limits_{\substack{ i=1 \\ i\neq k}}^{p+q}
\frac{\mathcal{B} \left( s, w_{i} \right) }{w_{i}^{2}} w_{i}.
\]
Using (\ref{cd02}) to get
\begin{eqnarray*}
  \varphi _{s} \left( w_{k} \right)  &=& -\frac{1}{s^{2}} \left( \left(
      \frac{\left( \mathcal{B} \left( s,w_{k} \right) \right)
        ^{2}}{w_{k}^{2}} - y_{k}^{2} \right) w_{k} + \left( 2
      \mathcal{B} \left( s,w_{k}\right) \right) y_{k}\right) , \\
  &=-&\frac{1}{s^{2}} \left( \left( \frac{\left( \mathcal{B} \left( s,
            w_{k} \right) \right)^{2}}{w_{k}^{2}} - \left( \sum\limits _{\substack{ i=1 \\ i\neq
            k}}^{p+q} \frac{ \left( \mathcal{B} \left( s, w_{i} \right) \right)
          ^{2}}{w_{i}^{2}} \right) \right) w_{k} + 2\mathcal{B} \left( s, w_{k}
    \right) y_{k} \right).
\end{eqnarray*}
From this result we obtain that the $k$-th column of the matrix $A_{s}
= \left[ \varphi _{s} \right] _{B}$ is given by

\begin{equation*}
\left[ \varphi _{s} \left( w_{k} \right) \right] _{B} =
-\frac{1}{s^{2}} \left(
\begin{array}{cc}
  2 \mathcal{B} \left( s, w_{k} \right) \frac{\mathcal{B} \left( s,
      w_{1} \right) }{w_{1}^{2}} &  \\
  2 \mathcal{B} \left( s, w_{k} \right) \frac{\mathcal{B} \left( s,
      w_{2} \right) }{w_{2}^{2}} &  \\
  \vdots  &  \\
  \frac{\left( \mathcal{B} \left( s, w_{k} \right) \right)
    ^{2}}{w_{k}^{2}} -\left( \sum\limits_{\substack{ i=1 \\ i\neq
        k}}^{p+q} \frac{\left( {\mathcal{B}} \left( s, w_{i} \right)
      \right) ^{2}}{w_{i}^{2}} \right)  & \leftarrow \;\; (k\textrm{-th
    row}) \\
  2\mathcal{B} \left( s, w_{k} \right) \frac{\mathcal{B} \left( s,
      w_{k+1} \right) }{w_{k+1}^{2}} &  \\
  \vdots  &  \\
  2\mathcal{B} \left( s, w_{k} \right) \frac{\mathcal{B} \left( s,
      w_{p+q} \right) }{w_{p+q}^{2}} &
\end{array}
\right) .
\end{equation*}
Equivalently
\begin{equation}
 \label{cd04}
\left( A_{s}\right) _{lj} = \left\{
\begin{array}{c}
  -2 \mathcal{B} \left( s, w_{j} \right) \frac{\mathcal{B} \left( s,
      w_{l} \right) }{s^{2} w_{l}^{2}} \textrm{ if } l\neq j , \\
  -\frac{1}{s^{2}} \left( \frac{\left( \mathcal{B} \left( s, w_{j} \right) \right)
      ^{2}}{w_{j}^{2}} - \left( \sum\limits_{\substack{ i=1 \\ i\neq
          j}}^{p+q} \frac{\left( \mathcal{B} \left( s,w_{i} \right) \right)
        ^{2}}{w_{i}^{2}} \right) \right) \textrm{ for } l=j.
\end{array}
\right. 
\end{equation}

\section{Examples}
\label{sec:examples}

Let $T$ be the orthogonal transformation on the space $\mathbb{R}^{2,3}$
represented by the matrix
\begin{equation*}
T_E= \left[
\begin{matrix}
  1 & 5 & 4 & 3 & 0 \\
  -5 & 1 & 3 & -4 & 0 \\
  4 & 3 & 1 & 5 & 0 \\
  3 & -4 & -5 & 1 & 0 \\
  0 & 0 & 0 & 0 & -1
\end{matrix}
\right]
\end{equation*}
with respect to the canonical basis $E = \left\{ e_{1}, e_{2}, e_{3},
  e_{4}, e_{5} \right\}$ of $\mathbb{R}^{2,3}$. Then we have $T \left(
  e_{i} \right) -e_{i} \neq 0$, \ $\left( T \left( e_{i} \right) -
  e_{i} \right) ^{2} = 0$ for $i=1, 2, 3, 4$, and $T \left( e_{5}
\right) - e_{5} = -2e_{5}$. We can take $c_{1}=e_{5}$.

It is easy to see that $\varphi _{c_{1}} T$ restricted to $\langle
e_{1}, e_{2}, e_{3}, e_{4} \rangle$ satisfies the conditions of Lemma
2.9. By Lemma 2.8 and the proof of Lemma 3.1 we obtain $c_{2} =
\varphi_{c_{1}} T \left( e_{1} \right) + e_{1} = 2e_{1} - 5e_{2} +
4e_{3} + 3 e_{4}$ and $c_{3} = e_{1}$. Then we have $\varphi _{c_{3}}
\varphi _{c_{2}} \varphi _{c_{1}} T \left( e_{i} \right) =e_{i}$, for
$i=1, 5$. Notice that $\varphi _{c_{2}} \varphi _{c_{1}} T$ does not
satisfy the hypothesis of Lemma 2.9 on $ \langle e_{1}, e_{2}, e_{3},
e_{4} \rangle$.

Now we have $\left( \varphi _{c_{3}} \varphi _{c_{2}} \varphi _{c_{1}}
  T \left( e_{i} \right) - e_{i} \right) ^{2} \neq 0$ for $i= 2, 3,
4$. We can take
\begin{equation*}
  c_{4} = \varphi _{c3} \varphi _{c_{2}} \varphi _{c_{1}} T \left(
    e_{2} \right) - e_{2} = \frac{25}{2} e_{2} - 7 e_{3} -
  \frac{23}{2}e_{4}. 
\end{equation*}
Then we have $\varphi _{c_{4}} \varphi _{c_{3}} \varphi _{c_{2}}
\varphi _{c_{1}} T \left( e_{i} \right) =e_{i}$ for $i=1, 2, 5$, and

$\left( \varphi _{c_{4}} \varphi _{c_{3}} \varphi _{c_{2}} \varphi
  _{c_{1}} T \left( e_{i} \right) - e_{i} \right) ^{2} \neq 0$ for
$i=3, 4$.

We can take $c_{5} = \varphi _{c_{4}} \varphi _{c3} \varphi _{c_{2}}
\varphi _{c_{1}} T \left( e_{3} \right) - e_{3} = \frac{-18}{25} e_{3}
+ \frac{24}{25}e_{4}$.

It is easy to verify that $\varphi _{c_{5}} \varphi _{c_{4}}
\varphi_{c_{3}} \varphi_{c_{2}} \varphi _{c_{1}} T \left( e_{i}
\right) = e_{i}$ for $i=1, 2, 3, 4, 5$, that is $T = \varphi _{c_{1}}
\varphi _{c_{2}} \varphi _{c_{3}} \varphi_{c_{4}} \varphi _{c_{5}}$.

Let us denote by $A_{j}$ the matrix representation with respect to the
canonical basis $E$ of $\varphi _{c_{j}}$, for $1 \le j \le 5$. 
Using formula (\ref{cd04}) we obtain:

\begin{equation*}
A_{5}=\left[
\begin{matrix}
  1 & 0 & 0 & 0 & 0 \\
  0 & 1 & 0 & 0 & 0 \\
  0 & 0 & \frac{7}{25} & \frac{24}{25} & 0 \\
  0 & 0 & \frac{24}{25} & -\frac{7}{25} & 0 \\
  0 & 0 & 0 & 0 & 1
\end{matrix}
\right], \quad A_{4}=\left[
\begin{matrix}
  1 & 0 & 0 & 0 & 0 \\
  0 & \frac{27}{2} & 7 & \frac{23}{2} & 0 \\
  0 & -7 & -\frac{73}{25} & -\frac{161}{25} & 0 \\
  0 & -\frac{23}{2} & -\frac{161}{25} & -\frac{479}{50} & 0 \\
  0 & 0 & 0 & 0 & 1
\end{matrix}
\right], \quad A_{3}=\left[
\begin{matrix}
  -1 & 0 & 0 & 0 & 0 \\
  0 & 1 & 0 & 0 & 0 \\
  0 & 0 & 1 & 0 & 0 \\
  0 & 0 & 0 & 1 & 0 \\
  0 & 0 & 0 & 0 & 1
\end{matrix}
\right],
\end{equation*}

\begin{equation*}
  A_{2}=\left[
\begin{matrix}
  -1 & 5 & 4 & 3 & 0 \\
  5 & -\frac{23}{2} & -10 & -\frac{15}{2} & 0 \\
  -4 & 10 & 9 & 6 & 0 \\
  -3 & \frac{15}{2} & 6 & \frac{11}{2} & 0 \\
  0 & 0 & 0 & 0 & 1
\end{matrix}
\right],  \quad A_{1}=\left[
\begin{matrix}
  1 & 0 & 0 & 0 & 0 \\
  0 & 1 & 0 & 0 & 0 \\
  0 & 0 & 1 & 0 & 0 \\
  0 & 0 & 0 & 1 & 0 \\
  0 & 0 & 0 & 0 & -1
\end{matrix}
\right].
\end{equation*}

A direct computation yields $T_E =A_{1} A_{2} A_{3} A_{4} A_{5}$.

Computing the Clifford product $c_1 c_2 c_3 c_4 c_5$ we obtain a
linear combination of $s$-vectors where $s \le 3$. This suggests that
$c_1 c_2 c_3 c_4 c_5$ may be equal to the Clifford product of 3
vectors and that $T$ could be factored as the product of 3
reflections.  We present next another way to factor $T$ that confirms
this conjecture.
  
We apply the factorization algorithm to $T$, but now using the
orthogonal basis $W = \left\{ w_{1}, w_{2}, w_{3}, w_{4}, w_{5}
\right\}$, where

$ w_{1} = e_{3} + e_{4} - e_{5}, w_{2} = e_{1} + e_{2}, w_{3} = e_{1}
+ e_{4} + 2e_{5}, w_{4}= e_{3} - e_{4}$ and $w_{5} = e_{1} - e_{2}$.

We obtain $\left( T\left( w_{i} \right) -w_{i} \right) ^{2} \neq 0$,
for $ 1 \le i \le 5$. Therefore we can take
\begin{eqnarray*}
  d_{1} &=&T\left( w_{1}\right) -w_{1}, \\
  \varphi _{d_{1}} T \left( w_{i} \right) -w_{i} &\neq &0\text{\ and\  }
  \left( \varphi_{d_{1}} T \left( w_{i} \right) -w_{i} \right) ^{2} \neq
  0, \text{\ for\  } i= 2, 3, 4, 5,
  \\
  d_{2} &=&\varphi _{d_{1}}T\left( w_{2}\right) -w_{2},\text{ } \\
  \varphi _{d_{2}} \varphi _{d_{1}} T \left( w_{i} \right) -w_{i}  &
  \neq & 0\text{\ and\  } \left( \varphi _{d_{2}} \varphi _{d_{1}} T
    \left( w_{i} \right) - w_{i} \right)^{2} \neq 0, \text{\ for\  }
  i= 3, 4, 5, \\
  d_{3} &=&\varphi _{d_{2}} \varphi _{d_{1}} T \left( w_{3} \right) -
  w_{3}, \text{ }
  \\
  0 &=&\varphi _{d_{3}} \varphi _{d_{2}} \varphi _{d_{1}} T \left(
    w_{i} \right) - w_{i}, \text{\ for\  }i=4,5.
\end{eqnarray*}
and therefore $T = \varphi _{d_{1}} \varphi _{d_{2}} \varphi
_{d_{3}}$.

Let $B_{i}$ denote the matrix representation of $\varphi _{d_{i}}$
with respect to the canonical basis $E$, for $i = 1, 2, 3$.  We have

$$B_{1}=\left[
\begin{matrix}
  \frac{51}{2} & -\frac{7}{2} & -\frac{35}{2} & \frac{35}{2} & -7 \\
  -\frac{7}{2} & \frac{3}{2} & \frac{5}{2} & -\frac{5}{2} & 1 \\
  \frac{35}{2} & -\frac{5}{2} & -\frac{23}{2} & \frac{25}{2} & -5 \\
  -\frac{35}{2} & \frac{5}{2} & \frac{25}{2} & -\frac{23}{2} & 5 \\
  7 & -1 & -5 & 5 & -1
\end{matrix}
\right],\qquad
 d_{1}=7e_{1}-e_{2}+5e_{3}-5e_{4}+2e_{5},$$

$$B_{2}=\left[
\begin{matrix}
  \frac{347}{9} & -\frac{104}{9} & -\frac{286}{9} & \frac{208}{9} &
  -\frac{26}{3} \\
  -\frac{104}{9} & \frac{41}{9} & \frac{88}{9} & -\frac{64}{9} &
  \frac{8}{3}
  \\
  \frac{286}{9} & -\frac{88}{9} & -\frac{233}{9} & \frac{176}{9} &
  -\frac{22}{3} \\
  -\frac{208}{9} & \frac{64}{9} & \frac{176}{9} & -\frac{119}{9} &
  \frac{16}{3}
  \\
  \frac{26}{3} & -\frac{8}{3} & -\frac{22}{3} & \frac{16}{3} & -1
\end{matrix}
\right] ,\qquad d_{2} = 26e_{1} - 8e_{2} + 22e_{3} - 16e_{4} +
6e_{5},$$

$$B_{3}=\left[
\begin{matrix}
  \frac{43}{18} & -\frac{25}{18} & -\frac{35}{18} & \frac{5}{18} &
  -\frac{5}{3}
  \\
  -\frac{25}{18} & \frac{43}{18} & \frac{35}{18} & -\frac{5}{18} &
  \frac{5}{3}
  \\
  \frac{35}{18} & -\frac{35}{18} & -\frac{31}{18} & \frac{7}{18} &
  -\frac{7}{3}
  \\
  -\frac{5}{18} & \frac{5}{18} & \frac{7}{18} & \frac{17}{18} & \frac{1}{3} \\
  \frac{5}{3} & -\frac{5}{3} & -\frac{7}{3} & \frac{1}{3} & -1
\end{matrix}
\right] ,\qquad d_{3} = -5 e_{1} + 5 e_{2} - 7 e_{3} + e_{4} - 6
e_{5}.$$

It is easy to verify that $T_E= B_{1} B_{2} B_{3}$, and notice that
$d_{1} d_{2} d_{3} = 6 c_{1} c_{2} c_{3} c_{4} c_{5}$.

The matrix representation of $T$ with respect to the basis $W$ is
\[
T_W = \left[
\begin{array}{rrrrr}
  \frac{1}{3} & 2 &\frac{4}{3} & -\frac{10}{3} & \frac{8}{3} \\
  3 & 1 & 3 & 4 & -5\\
    \frac{2}{3} & 1 & -\frac{1}{3} & -\frac{5}{3} & \frac{4}{3}\\
    5 & 4 & 5 & 1 & -3\\
    4 & 5 & 4 & -3 & 1
\end{array}
\right].
\]

Note that this second factorization avoids the Artinian case, that is,
the situation where the hypothesis of Lemma 2.9 holds.

As another application of formula (\ref{cd04}), we can find the matrix
representation of each $\varphi_{d_{i}}$, for $i = 1, 2, 3$, but now
with respect to the orthogonal basis $W$. Using the above values of
$d_1$, $d_2$ and $d_3$, respectively, we obtain:
\[
C_{1}=\left[
\begin{array}{ccccc}
  \frac{1}{3} & -2 & \frac{4}{3} & \frac{10}{3} & -\frac{8}{3} \\
  3 & 10 & -6 & -15 & 12 \\
  \frac{2}{3} & 2 & -\frac{1}{3} & -\frac{10}{3} & \frac{8}{3} \\
  5 & 15 & -10 & -24 & 20 \\
  4 & 12 & -8 & -20 & 17
\end{array}
\right] ,
\]
\[
C_{2}=\left[
\begin{array}{ccccc}
  1 & 0 & 0 & 0 & 0 \\
  0 & 10 & -9 & -19 & 17 \\
  0 & 3 & -2 & -\frac{19}{3} & \frac{17}{3} \\
  0 & 19 & -19 & -\frac{352}{9} & \frac{323}{9} \\
  0 & 17 & -17 & -\frac{323}{9} & \frac{298}{9}%
\end{array}
\right] ,
\]
and
\[
C_{3}=\left[
\begin{array}{ccccc}
  1 & 0 & 0 & 0 & 0 \\
  0 & 1 & 0 & 0 & 0 \\
  0 & 0 & -2 & -\frac{4}{3} & \frac{5}{3} \\
  0 & 0 & -4 & -\frac{7}{9} & \frac{20}{9} \\
  0 & 0 & -5 & -\frac{20}{9} & \frac{34}{9}
\end{array}
\right].
\]

A direct computation yields $T_W = C_1 C_2 C_3$. It should be noticed
that the orthogonal basis $W$ is not ordered according to the
signature of $\mathbb{R}^{2,3}$.

\section{Final remarks and conclusions}

The factorization of an orthogonal transformation is not unique, as we
have shown in our example. Actually, the number of reflections that
factorizes an orthogonal transformation is not unique either. An
interesting question related to this last point is to determine the
minimum number of reflections by hyperplanes required to factorize a
given orthogonal transformation. An answer to this question using
matrices has been given in \cite{S}, and it was translated to the
language of orthogonal transformations in \cite[Pgs 260,261]{ST}.

In our examples (Section \ref{sec:examples}) we use the grade of a
multivector to find the minimum number of reflections by hyperplanes
required to factorize $T_E$. We can formalize this procedure as
follows:

\begin{lem}
  \label{lem:minimum}
  Let $T$ be an orthogonal transformation on
  $\mathbb{R}^{p,q}$. Assume that there exist invertible elements
  $s_1, \ldots, s_k$ in $\mathbb{R}^{p,q}$, with $k \le p+q$, such
  that
\begin{equation*}
  T \left( x \right) =\left( -1 \right) ^{k} s_{1} s_{2} \cdots
  s_{k-1} s_{k} x s_{k}^{-1} s_{k-1}^{-1} \cdots  s_{2}^{-1} s_{1}^{-1}, \qquad x
  \in \mathbb{R}^{p,q}.
\end{equation*}
If $gr\left(s_{1} \cdots s_{k-1} s_{k} \right)=t \leq k$, then $T$
cannot be factored into less than $t$ reflections with respect to
hyperplanes.
\end{lem}

\textit{Proof:} Suppose that there exist $s_1^\prime, \ldots,
s_l^\prime$, with $l<t$, such that
\begin{equation*}
  T \left( x \right) = \left( -1 \right) ^{k} s_{1}^\prime \cdots
  s_{l-1}^\prime s_{l}^\prime x ( s_{l}^\prime )^{-1}
  ( s_{l-1}^\prime )^{-1} \cdots ( s_1^\prime )^{-1}, \qquad x \in
  \mathbb{R}^{p,q}.
\end{equation*}
We have that $\textrm{gr} \left(s_{1}^\prime \cdots s_{l-1}^\prime
  s_{l}^\prime \right) \leq l<t$. But $\left(s_{1}^\prime \cdots
  s_{l-1}^\prime s_{l}^\prime \right)= \lambda \left (s_{1} \cdots
  s_{k-1} s_{k} \right)$, where $\lambda \in \mathbb{R}$ and therefore
$\textrm{gr} \left( s_{1}^\prime \cdots s_{l-1}^\prime s_{l}^\prime
\right) = \textrm{gr} \left(s_{1} \cdots s_{k-1} s_{k} \right)=t$,
which is a contradiction.  \mbox{$\blacksquare$}

With this result, we can state the following theorem

\begin{thm}
  Let $T$ be an orthogonal transformation on the space
  $\mathbb{R}^{p,q}$. If there exist invertible elements $s_1, s_2,
  \ldots, s_k \in \mathbb{R}^{p,q}$, such that:
 \begin{itemize}
 \item $ T\left( x \right) = \left( -1 \right) ^{k} s_{1} \cdots
   s_{k-1} s_{k} x s_{k}^{-1} s_{k-1}^{-1} \cdots s_{2}^{-1}
   s_{1}^{-1}, \qquad x \in \mathbb{R}^{p,q}$,
 \item $\textrm{gr} \left( s_{1}, \dots, s_{k-1}, s_{k} \right) = t$,
   and
  \item $Ker(T-I)$ is non-degenerate,
 \end{itemize}
 then $T$ can be factored into $t$ reflections with respect to
 hyperplanes and, moreover, $\textrm{dim} \left( \textrm{Ker} \left( T
     - I \right) \right) ^{\perp}= t$.
\end{thm}
\textit{Proof:} By hypothesis $T = \varphi_{s_1} \cdots \varphi_{s_k}$
and $\textrm{gr} \left( s_{1}, \dots, s_{k-1}, s_{k} \right) = t$.  By
considering the space $V_1 = \langle s_{1}, \cdots, s_{k} \rangle$,
where $\dim V_1=t$, it is easy to show that $V_1^{\bot} \subset
\textrm{Ker} \left( T - I \right)$. If $\textrm{Ker} \left( T - I
\right)$ is non-degenerate. Then we can find an orthogonal basis
$\mathcal{B} = \left\{w_1, \dots, w_j, w_{j+1} , \dots,
  w_{j+l}\right\}$, where $j+l=n$ and $\textrm{Ker} \left( T-I
\right)=\langle w_{1}, \dots, w_{j-1}, w_{j} \rangle $ and applying
the algorithm proposed in this paper we get $\varphi_i = I$, for $i
=1, 2, \dots, j$. From the application of the \CD theorem to $T
|_{V_{j+1}}$, where $V_{j+1} = \langle w_{j+1}, \dots, w_{j+l}
\rangle$, we obtain that $T$ is the composition of at most $l$
reflections by hyperplanes. From Lemma 6.1 we know that $t \leq l$.

Now, suppose that there exist $u_1, \ldots, u_m \in \left(\textrm{Ker}
  \left(T-I\right) \right)^{\bot}$, such that $T = \varphi_{u_1}
\cdots \varphi_{u_m}$. If $ \textrm{gr} \left(u_1 \cdots u_m \right) <
l$ then by considering $T |_{V_{j+1}}$ we can find a nonzero vector $u
\in V_{j+1}$ orthogonal to the set $\left \{u_1, \ldots u_m \right
\}$. But in such case $u \in \textrm{Ker} \left( T - I \right) \cap
\left(\textrm{Ker} \left( T - I \right) \right)^\bot$, which is a
contradiction since $\textrm{Ker} \left(T-I\right)$ is non-degenerate.
Therefore $T$ is the composition of $l$ reflections through
hyperplanes, that is, there exist $u_1, \ldots, u_l$ such that $T =
\varphi_{u_1} \cdots \varphi_{u_l}$ and $\textrm{gr} \left(u_1 \cdots
  u_l \right)=l$ and, moreover, $t=gr \left(s_1s_2 \cdots s_k \right)=
gr \left(u_1 u_2 \cdots u_l \right)=l= \dim \left(\textrm{Ker}
  \left(T-I\right) \right)^\bot$ \mbox{$\blacksquare$}

It should be remarked that the previous Theorem relates the grade of a
multivector with the Cartan-Diedudonn\'e-Scherk theorem. It remains
however, to propose an algorithm capable of finding explicitly the
minimum number of reflections required to decompose a given orthogonal
transformation.

\begin{ack}
  MAR would like to thank for financial support from COFAA-IPN.  GAG
  was supported by the Program for the Professional Development in
  Automation, through the grant from the Universidad Aut\'{o}noma
  Metropolitana and Parker Haniffin - M\'{e}xico. JLA would like to
  thank DGAPA-UNAM (grant IN100310-3) and CONACyT (grant 50368 and
  79641) for financial support.
\end{ack}

\end{document}